\documentclass[preprint,12pt]{elsarticle}

\usepackage{graphicx}

\usepackage{amsthm}
\usepackage{amsmath}
\usepackage{amssymb}

\journal{Journal of Computational and Applied Mathematics}

\def\R{\mathbb{R}}
\newcommand{\bh}{$h$}
\newcommand{\bp}{$p$}
\newcommand{\bhp}{$hp$}

\begin{document}

\begin{frontmatter}



\title{An \bhp-adaptive strategy for elliptic problems}


\author{Hui Liu, Tao Cui, Wei Leng, Linbo Zhang}

\address{ State Key Laboratory of Scientific and Engineering Computing \\
Academy of Mathematics and Systems Science \\
Chinese Academy of Sciences \\
P.O. Box 2719, Beijing 100190, CHINA}

\begin{abstract}
In this paper a new \bhp-adaptive strategy for elliptic problems
based on refinement history is proposed, which chooses
\bh-, \bp- or \bhp-refinement on individual elements
according to a posteriori error estimate, as well as
smoothness estimate of the solution obtained by comparing the
actual and expected error reduction rate.
Numerical experiments show that exponential convergence
can be achieved with this strategy.
\end{abstract}


\end{frontmatter}


\section{Introduction}

The adaptive finite element method (AFEM) is a widely used
numerical method for solving partial differential equations.
The \bh-version of AFEM modifies the size of the elements
(\bh-refinement) while keeping the polynomial degrees fixed
\cite{houston}. The \bp-version of AFEM adjusts the polynomial
degrees in the elements (\bp-refinement) while keeping the
size of the elements fixed. The \bhp-version of AFEM is more general,
which consists of combining freely \bh-refinement and \bp-refinement.
The \bhp-version of AFEM  dates back to 1986, thanks to the pioneering
work of Ivo Babu\v{s}ka et al. \cite{suri,guo1,babu1,babu2,gui}.
With \bhp~AFEM exponential convergence could be achieved if \bh-refinement
and \bp-refinement are integrated properly~\cite{melenk,anis,guo1}.

One essential issue in the \bhp-adaptive finite element method is
the design of refinement strategy, i.e., to decide which element
should be refined and which kind of refinement should be performed.
According to approximation theory,
\bp-refinement should be performed on elements in which the
solution to the partial differential equations is smooth and
\bh-refinement should be performed on elements in which the
solution is non-smooth~\cite{melenk}. Unfortunately, since the
property of the solution is usually unknown, we need to
estimate its smoothness using the computed numerical solution and
other data. For this purpose many strategies have been proposed and
developed. Owens et al.~\cite{bernardi,valen} used a priori
information of the computational domain and boundary data
to determine the location of singularities of the solution,
and performed \bh-refinement on elements which had
singularities and \bp-refinement elsewhere.
Oden et al.~\cite{oden} introduced the so-called Texas-3-step strategy.
Melenk et al.~\cite{melenk} and Heuveline et al.~\cite{rannacher} proposed
heuristic strategies which made use of the refinement history.
Another class of strategies consisted of using
error estimators obtained from solving local problems as
indicators for guiding the refinement~\cite{anis,dem2}.
For other strategies proposed and studied in the literature,
we refer to~\cite{affia,gui,melenk1,anis,paul,mavr}.

In this paper, we propose an \bhp-refinement strategy which is based
on a posteriori error estimate and estimation of the smoothness of the
solution using the reduction rates of the a posteriori error estimate
in the refinement history.
%
%
This strategy is mainly motivated by Melenk et al.~\cite{melenk} and Heuveline
et al.~\cite{rannacher}, it removes the requirement of regular refinement
and the dependence on mesh size $h$ in~\cite{melenk,rannacher},
and can be applied to both two and three dimensional elliptic problems.

The layout of the paper is as follows. In \S 2, the model problem
and notations are introduced. In \S 3, the \bhp-adaptive strategy
is deduced in details. In \S 4, the efficiency of the new strategy
is illustrated and compared to some other strategies through two
numerical examples. In \S 5, some concluding remarks are given.

\section{Model problem and notations}

For a bounded Lipschitz domain $\Omega \subset \R^d$, $d = 2, 3$, the
following model problem is considered:
\begin{equation}
\label{poisson}
 -\Delta u = f \text{ on } \Omega, \quad u = g\text{ in } \partial \Omega.
\end{equation}
where $f \in L^2(\Omega)$. The problem can be read in the weak form:
find $u \in H_0^1(\Omega)$ such that
\begin{equation}
 a(u, v) = L(v), \quad \forall v \in H_0^1(\Omega),
\end{equation}
where
\begin{equation}
a(u, v) = \int_{\Omega}\nabla u\cdot \nabla v dx, \quad
L(v) = \int_{\Omega} fvdx .
\end{equation}

Our goal is to design an \bhp-finite element subspace
$V_{hp} \subset H_0^1(\Omega)$ and to compute a numerical solution
$u_{hp} \in V_{hp} $ such that
\begin{equation}
a(u_{hp}, v_{hp}) = L(v_{hp}), \quad \forall v_{hp} \in V_{hp},
\end{equation}
and the error meets prescribed tolerance. Here for simplicity of
description we will assume $g=0$. In this case $g\ne0$, the problem
can be easily converted to the case $g=0$ with a shift operator.

For the sake of convenience, some notations are introduced here.
In the subsequent descriptions, we will denote by $u$ the exact solution
of Problem~(\ref{poisson}), by $u_{hp}$ the numerical solution of the problem
with respect to a triangulation $T$ and a finite element space $V_{hp}$ on $T$,
and by $e=u - u_{hp}$ the error between the exact solution and the
numerical solution. The energy norm,
$\|\cdot\|$, is defined as $\| u \| = \sqrt{a(u, u)}$. In the
adaptive process, $\varepsilon$ stands for the tolerance, which is the
stop criterion, $\eta_K$ stands for the error indicator defined on element $K$,
and $\eta = (\sum_{K \in T}\eta_K^2)^{1/2}$ is the global error indicator.
For a given element $K$, $p_K$ and $h_K$ denote the degree of the polynomial
basis functions on $K$ and the diameter of $K$, respectively.
When the element $K$ is divided (refined) into subelements,
$c_K$ denotes the number of its children. Finally, $N_d$ is used to
denote the total number of degrees of freedom in the mesh $T$.

\section{An \bhp-adaptive strategy}
In this section we give our \bhp-adaptive strategy. This strategy is
based on the \emph{expected error reduction factors} of \bh-, \bp-,
or \bhp-refinement.
The expected error reduction factors are calculated under the assumption
that the numerical solution converges algebraically under \bh-refinement
and exponentially under \bp-refinement. We will first deduce the expected
error reduction factors for various refinement types, then describe the
new \bhp-adaptive strategy in details.

First we deduce the expected error reduction factor $\lambda_h$ for
\bh-refinement. We assume that the optimal convergence rate of the
\bh-version of adaptive finite element method is algebraic, which
can be written as \cite{zhoua,fernan,cascon},
\begin{equation}
\label{converge1}
 \|e\| \le C_1N_d^{-\frac{p}{d}},
\end{equation}
where $p$ is the degree of the piecewise polynomials. Suppose the
fine mesh $T_1$ is obtained from uniform refinement of the mesh
$T$ by dividing each element $K$ into $c_K$ subelements. Then the
number of degrees of freedom on mesh $T_1$ is about $c_KN_d$.

Suppose we have an appropriate error indicator $\{\eta_K \mid K\in T\}$.
We make the following hypotheses.

\bigskip\noindent
{\bf (H1)} The error indicator is precise, i.e., there exist
constants $C_1$ and $C_2$ such that,
\begin{equation}
\label{converge2} \|e\| = C_1N_d^{-\frac{p}{d}} = C_2(\sum_{K \in
T}\eta_K^2)^{\frac{1}{2}},
\end{equation}
{\bf (H2)}  For any element $K$, the error indicators on all its children
are equal.
\bigskip

Let $\lambda_h$ be the expected error reduction factor for \bh-refinement.
By combining (H1) and (H2), we get the following relationship
\begin{equation}
\label{converge3} C_1(c_KN_d)^{\frac{-p}{d}} =C_2(\sum_{K' \in
T_1}\eta_{K'}^2)^{\frac{1}{2}} =C_2(c_K\lambda_h^2\sum_{K \in
T}\eta_K^2)^{\frac{1}{2}}.
\end{equation}
Comparing ~(\ref{converge3}) to
~(\ref{converge2}), we have
\begin{equation}
\lambda_h^2 = \frac{1}{c_K}(\frac{1}{c_K})^{\frac{2p}{d}}.
\end{equation}
To improve the efficiency, we use a slightly enlarged $\lambda_h$,
which is given by
\begin{equation}
\label{lbdh} \lambda_h = (\frac{1}{c_K})^{\frac{p}{d}}.
\end{equation}

Next we deduce the expected error reduction factor $\lambda_p$
for \bp-refinement. In \bp-refinement the mesh is fixed and the degree
of the polynomials is adjusted. On a quasi-uniform mesh with uniform
polynomial degree the following error estimation is expected~\cite{anis,suri}
\begin{equation}
\label{pcon}
\|e\|_{H^1(\Omega)} \le C\frac{h^\mu}{p^{(m-1)}} \|u\|_{H^m(\Omega)},
\end{equation}
where $h$ is the mesh size, $p$ the polynomial degree, $\mu = \min
(p, m-1)$, $C$ a constant independent of $h$ and $p$, and $u
\in H^m(\Omega)$. We make the following hypothesis.

\bigskip
\noindent
{\bf (H3)}
$\|e\|_{H^1(\Omega)} = C\frac{h^\mu}{p^{(m-1)}} \|u\|_{H^m(\Omega)}$
and $p \ge (m-1)$.
\bigskip

When the degree $p$ is increased by one, by (H3) we have
\begin{equation}
 \|e\|_{H^1(\Omega)} =  C\frac{h^\mu}{(p+1)^{(m-1)}} \|u\|_{H^m(\Omega)}
 = C(\frac{p}{p+1})^{m-1}\frac{h^\mu}{p^{(m-1)}} \|u\|_{H^m(\Omega)}.
\end{equation}
Thus the error reduction factor $\lambda_p$ is
\begin{equation}
\lambda_p = (\frac{p}{p+1})^{m-1}.
\end{equation}
$m$ is a positive integer satisfying (H3). In this
paper we set $m$ to $p/2+1$. Then we have
\begin{equation}
\label{lbdp} \lambda_p = (\frac{p}{p+1})^{\frac{p}{2}}.
\end{equation}

Finally the expected error reduction factor $\lambda_{hp}$ for
\bhp-refinement can readily be obtained by combining $\lambda_h$ and
$\lambda_p$, which is given by
\begin{equation}
\label{lbdhp} \lambda_{hp} = (\frac{p}{p+1})^{\frac{p}{2}}(\frac{1}{c_K})^{\frac{p}{d}}.
\end{equation}

As a widely accepted criterion in adaptive finite element
methods, the error should be distributed asymptotically uniformly over all
elements~\cite{melenk}. Therefore, elements with large error estimator should
be marked for refinement. Here we employ the so-called \emph{maximum strategy},
which can be described as follows
\begin{equation}
\eta_K \ge \alpha \max_{K' \in T}\eta_{K'}
\Leftrightarrow K \text{ is marked for refinement},
\end{equation}
where $\alpha \in (0, 1)$ is a predetermined parameter.

Our \bhp-adaptive strategy is given below, which is motivated
by Heuveline et al.~\cite{rannacher} and Melenk et al.~\cite{melenk}, using
a similar framework. Here \bh-refinement means dividing the element
into $c_K$ subelements, \bp-refinement means increasing the polynomial
degree by 1.

\medskip
\begin{enumerate}
\item[{\bf Step 1:}] Solve the problem on the current mesh $T$ with the current
setting of polynomial orders and compute the error indicator
$\{\eta_K\mid K\in T\}$ and the global error indicator $\eta$.
The adaptive process is stopped if $\eta$ is less than or equal to
$\varepsilon$ on the current mesh.

\item[{\bf Step 2:}] Mark elements for refinement using maximum strategy.

\item[{\bf Step 3:}] For each marked element $K$:
  \begin{itemize}
     \item If element $K$ is obtained by \bh-refinement of its parent
     element $K_m$, then check whether the following condition holds
     $$\eta_K^2 \le \lambda_h^2 \eta_{K_m}^2.$$
     If yes then mark $K$ for \bp-refinement.
     Otherwise mark $K$ for \bh-refinement.

     \item If element $K$ is obtained by \bp-refinement of its parent
     element $K_m$, then check whether the following condition holds
     $$\eta_K^2 \le \lambda_p^2 \eta_{K_m}^2.$$
     If yes then mark $K$ for \bp-refinement.
     Otherwise mark $K$ for \bh-refinement.

     \item If element $K$ is obtained by \bhp-refinement of its
     parent element $K_m$, then check whether the following condition holds
     $$\eta_K^2 \le \lambda_{hp}^2 \eta_{K_m}^2.$$
     If yes then mark $K$ for \bp-refinement.
     Otherwise mark $K$ for \bh-refinement.

     \item If element $K$ is not refined in the preceding adaptive step,
     then mark $K$ for \bp-refinement.
  \end{itemize}
\item[{\bf Step 4:}] Perform \bh-, \bp- or \bhp-refinement
as determined by {\bf Step 3}.\footnote{When we perform \bh-refinement
additional elements may be refined in order to maintain the conformity
of the mesh.}
\item[{\bf Step 5:}] Go to {\bf Step 1}.
\end{enumerate}

\medskip
The underlying idea behind the above process is that because of the exponential
convergence rate of \bp-refinement, it is preferred over \bh-refinement
whenever the solution is smooth. If the expected error reduction factor
is achieved in the previous refinement, then the solution is considered
smooth and \bp-refinement is performed, otherwise \bh-refinement is performed.

Remark: the strategy proposed by Melenk et al.~\cite{melenk} was designed for two
dimensional problems. Our strategy is suitable for both two and three
dimensional problems and different error
reduction factors are deduced. For the strategy proposed by Rannacher et al.
\cite{rannacher}, the error reduction factor depended on the size of
elements. This dependency is removed in this paper.

\section{Numerical results}
In this section two examples are employed to illustrate the
efficiency of the new \bhp-adaptive strategy. These examples are
also computed using a traditional \bh-version adaptive finite element
method and another existing \bhp-adaptive strategy for comparison.

We have implemented our new \bhp-adaptive strategy using the
parallel adaptive finite element toolbox PHG~\cite{phg}. The
computations were performed on the cluster LSSC-III of the
State Key Laboratory of Scientific and Engineering Computing,
Chinese Academy of Sciences.

In these examples, since bisection refinement is used for \bh-refinement,
we have $c_K=2$, thus the expected error reduction
factors are given by
\begin{equation}
\label{lambdah} \lambda_h = (\frac{1}{2})^{\frac13 p_K},
\end{equation}

\begin{equation}
\label{lambdap} \lambda_p = (\frac{p_K - 1}{p_K})^{\frac12(p_K - 1)},
\end{equation}

\begin{equation}
\label{lambdahp} \lambda_{hp} = (\frac{p_K - 1 }{p_K})^{\frac12(p_K -
1)}(\frac{1}{2})^{\frac13(p_K - 1)}.
\end{equation}

The error indicator used here is the one introduced by Melenk et al.
\cite{melenk}. Though it was designed for two dimensional problems, it
is also valid for three dimensional problems. This error indicator is given by
\begin{equation}
\label{error} \eta_K^2 = \frac{h_K^2}{p_K^2}\|f_{p_K} + \Delta
u_{hp}\|_{L^2(k)}^2 + \sum_{f \subset \partial
K\cap\Omega}{\frac{h_f}{2p_f} \|[\frac{\partial u_{hp}}{\partial
n_f}]\|^2_{L^2(f)}},
\end{equation}
where $f_{p_K}$ is the $L^2(K)$-projection of the function $f$ on the space of
polynomials of degree $p_K - 1$, $h_f$ denotes the diameter of the face $f$,
$p_f = \max (p_{K_1}, p_{K_2})$, where $K_1$ and
$K_2$ are the two elements sharing the face $f$, and $[\cdot]$ denotes
the jump of a function across the face $f$.

The parameter $\alpha$ in the maximum strategy is chosen as $0.5$.
The linear systems of equations are solved by the PCG (Preconditioned
Conjugate Gradient) method with a block Jacobi preconditioner.
The initial meshes are generated using NETGEN~\cite{netgen} and
the initial polynomial degrees on all elements are set to 2.

For three dimensional Poisson equation  the optimal convergence rate
is exponential and is expected to be~\cite{guo1}
\begin{equation}
\|e\| \le C\exp\bigl(-\gamma (N_d)^{1/5}\bigr),
\end{equation}
where $\gamma$ is a constant.

In the figures the logarithm of the energy
error is plotted against $(N_d)^{1/5}$, and
three different strategies are compared.
The first one is a traditional \bh-adaptive finite
element method, denoted by ``HAFEM''. The second one is the
\bhp-adaptive strategy introduced in this paper, denoted by ``HP/PHG''. The last
one is the strategy of Melenk et al., denoted by ``HP/MK''.

\textbf{Example 4.1.} In this example, the domain is an $L$-shaped domain
given by $\Omega = (-1, 1)^3 \setminus (0,1]\times[-1,0)\times(-1,1)$,
and the analytic solution is given by
$u = \cos(2 \pi x)\cos(2 \pi y)\cos(2 \pi z)$.
The main difficulty in applying high order finite element methods to
this problem is that the even and odd derivatives of the solution behave
differently at each point in the domain, hence pure \bp-refinement may
not improve the numerical solution~\cite{zib}. The initial mesh is
uniform with 144 elements.

\begin{figure}[htb]
\centering
\includegraphics[width=0.8\textwidth]{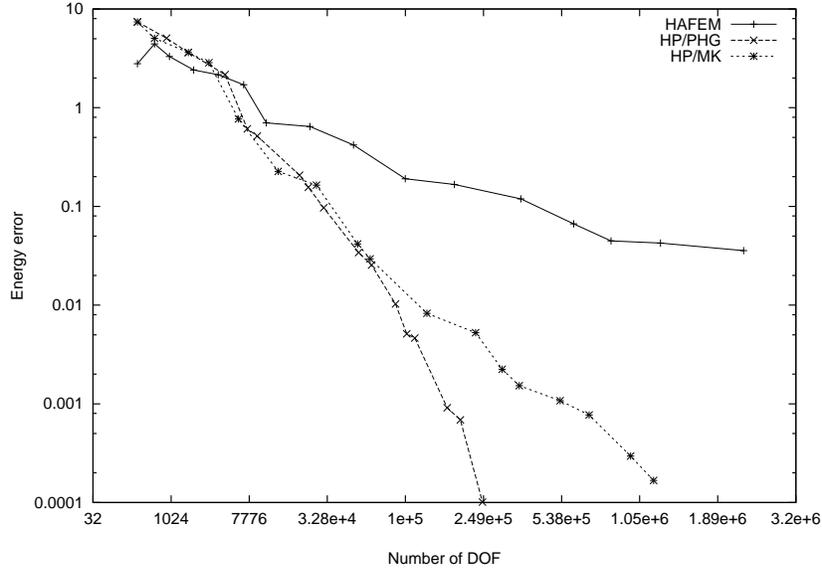}
\caption{Convergence histories (Example 4.1)}
\label{fex2}
\end{figure}

The convergence histories of different strategies are shown in Figure~\ref{fex2}
and statistics about the final meshes are shown in Table~\ref{tex2}.
We can observe that
the two \bhp~strategies exhibit exponential convergence rate while the
\bh-version converges algebraically. We can also observe that the
HP/PHG strategy performs better than the HP/MK strategy.

\begin{table}[htb]
\centering
\caption{The final meshes (Examples 4.1)}
\begin{tabular}{|l|r|r|c|} \hline
     & \# elements  & \# DOF   & Energy error   \\ \hline
HP/PHG & 3,772       & 246,046   &    1.01e-4    \\ \hline
HP/MK  & 35,696      & 1,171,216 &    1.67e-4     \\ \hline
HAFEM  & 1,663,068   & 2,263,137 &    3.57e-2     \\ \hline
\end{tabular}
\label{tex2}
\end{table}

\textbf{Example 4.2.} In this example,
the computational domain is given by $\Omega = (-1, 1)^3 \setminus [0, 1)^3$,
and the analytic solution is given by $u = {(x^2 + y^2 + z^2)}^{\frac{1}{4}}$,
whose gradient has a vertex singularity.
The initial mesh is uniform with 172 elements.

\begin{figure}[htb]
\centering
\includegraphics[width=0.8\textwidth]{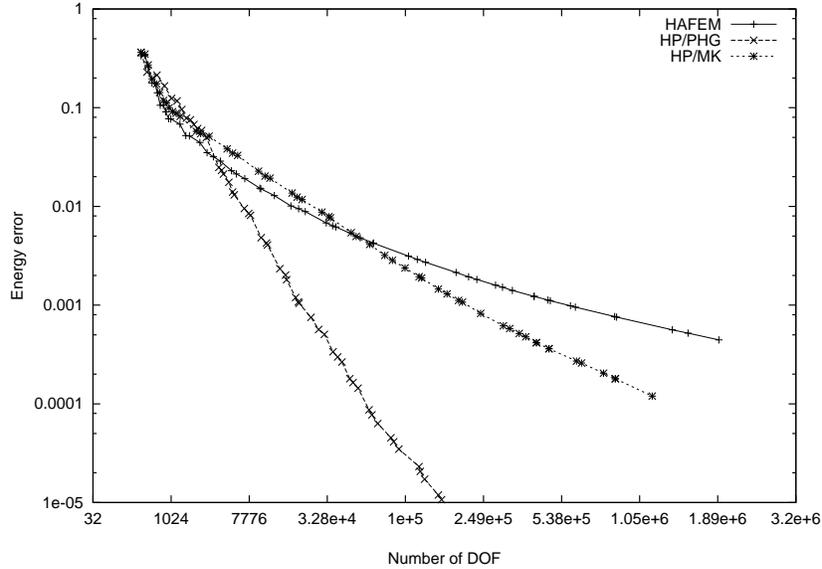}
\caption{Convergence histories (Examples 4.2)}
\label{fex3}
\end{figure}

The convergence histories and final meshes are shown in
Figure~\ref{fex3} and Table~\ref{tex3} respectively. Again for this
example, the \bh-version converges algebraically while the two
\bhp-versions converge exponentially. Data in Table~\ref{tex3} shows
that the performance of our strategy is much better than that of the
HP/MK strategy.

\begin{table}[htb]
\centering
\caption{Final meshes (Example 4.2)}
\begin{tabular}{|l|r|r|c|} \hline
     & \# elements  & \# DOF   & Energy error   \\ \hline
HP/PHG & 3,429       & 155,812   &  1.07e-5     \\ \hline
HP/MK  & 163,204     & 1,158,279 &  1.20e-4     \\ \hline
HAFEM  & 1,377,588   & 1,904,054 &  4.44e-4     \\ \hline
\end{tabular}
\label{tex3}
\end{table}

\section{Conclusion}

A simple and easy to implement \bhp-adaptive strategy based on error
reduction prediction is proposed. This strategy is suitable for two and
three dimensional problems. The efficiency of the strategy is demonstrated
through two numerical examples. Although the strategy is discussed
with the Poisson equation in this paper, it is applicable to general
elliptic problems. It also provides a general framework which can be easily
extended to other problems.

\section*{Acknowledgments}
This work is supported by the 973 Program under the grant
 2011CB309703, by China NSF under the grants 11021101 and 11171334, by the 973 Program under
 the grant 2011CB309701, the China NSF under the grants 11101417
and by the National Magnetic Confinement Fusion Science Program under the grants 2011GB105003.

\bibliographystyle{elsarticle-num}
\bibliography{<your-bib-database>}



\end{document}